\documentclass{article}

\usepackage{amssymb}
\usepackage{amsmath,amsthm}
\usepackage[dvips]{graphicx}
\usepackage{wrapfig}
\usepackage{bm}

\newtheorem{theorem}{Theorem}[section]
\newtheorem{corollary}[theorem]{Corollary}
\newtheorem{lemma}[theorem]{Lemma}
\newtheorem{proposition}[theorem]{Proposition}

\theoremstyle{definition}

\theoremstyle{remark}

\begin{document}
\title{Prime alternating knots of minimal warping degree two}

\author{Ayaka Shimizu 
\thanks{Department of Mathematics, National Institute of Technology, Gunma College, 580 Toriba-machi, Maebashi-shi, Gunma, 371-8530, Japan. Email: shimizu@gunma-ct.ac.jp }}

\maketitle

\begin{abstract}
The warping degree of an oriented knot diagram is the minimal number of crossing changes which are required to obtain a monotone diagram from the diagram. 
The minimal warping degree of a knot is the minimal value of the warping degree for all oriented minimal diagrams of the knot. 
In this paper, all prime alternating knots with minimal warping degree two are determined. 
\end{abstract}

\section{Introduction}

A {\it knot} is an embedding of a circle in $S^3$. 
A {\it knot projection} is a projection of a knot on $S^2$ such that each intersection is a double point where arcs cross transversely. 
We call such double point a {\it crossing}. 
A {\it knot diagram} is a knot projection with over/under information at each crossing. 
We denote by $|D|$ the knot projection which is obtained from a knot diagram $D$ by forgetting the over/under information. 
An {\it edge} of $D$ or $|D|$ is an arc from a crossing to the next crossing. 

For an oriented knot diagram $D$, take a base point $b$ on an edge. 
We call it a {\it based diagram} and denote it by $D_b$. 
A crossing $p$ of $D$ is a {\it warping crossing point} of $D_b$ if we encounter $p$ as an under-crossing first by traveling $D$ from $b$ with the given orientation. 
The {\it warping degree}, $d(D_b)$, of an oriented based diagram $D_b$ is the number of warping crossing points of $D_b$. 
The {\it warping degree}, $d(D)$, of an oriented diagram $D$ is the minimal value of $d(D_b)$ for all base points $b$. 
(\cite{kawauchi-lecture, shimizu-wd}. See also \cite{LM, MO} for similar notions.)
A knot diagram $D$ with $d(D)=0$ is said to be {\it monotone}. 
Every monotone diagram represents the trivial knot. 

The warping degree relates to some basic knot invariants. 
The {\it unknotting number}, $u(K)$, of a knot $K$ is the minimal number of crossing changes which are required to make $K$ unknotted. 
We have $u(K) \leq d(D)$ for any oriented diagram $D$ of $K$. 
The {\it ascending number}, $a(K)$, of a knot $K$ is defined in \cite{MO}, and is equivalent to the minimal value of $d(D)$ for all oriented knot diagrams $D$ of $K$.\footnote{
The relations of the ascending number and the unknotting number, crossing number and bridge number are shown in \cite{MO}. }
We also have $a(K) \leq d(D)$ for any oriented diagram $D$ of $K$ by definition. 
The strong point of the warping degree is that we can obtain the value easily by traveling a knot diagram.\footnote{
Applying the {\it warping-degree labeling}, we can obtain the value of the warping degree more easily. See \cite{shimizu-w-poly}. } 
In \cite{shimizu-wd}, the relation between the warping degree and the crossing number is given as follows. 
Let $D$ be an oriented diagram of a knot $K$, and let $-D$ denote $D$ with orientation reversed. 
Let $e(K)$ be the minimal value of $d(D)+d(-D)$ for all minimal diagrams $D$ of $K$. 
The inequality $e(K) \leq c(K)-1$ holds, and the equality holds if and only if $K$ is a prime alternating knot, where $c(K)$ is the crossing number of $K$. 
The knot invariant $e(K)$ is named the {\it warping sum} and studied in \cite{JS}. 
There are various breakdowns of the warping sum of $e(K)=d(D)+d(-D)$ for minimal diagrams $D$ of a prime alternating knot $K$ although the value is constant to $c(K)-1$. 
Then our next interest is how small the warping degree can be for minimal diagrams of a prime alternating knot. 

The {\it minimal warping degree}, $md(K)$, of a knot $K$ is the minimal value of $d(D)$ for all oriented minimal diagrams $D$ of $K$. 
By definition, $a(K) \leq md(K)$ holds for each knot $K$. 
In \cite{JS}, all the knots of minimal warping degree one are determined to be $3_1$ and $4_1$.\footnote{In \cite{MO}, it is shown that the ascending number of a knot $K$ is one if and only if $K$ is a twist knot. } 
In this paper, all the prime alternating knots of minimal warping degree two are determined as follows: 

\phantom{x}
\begin{theorem}
All the prime alternating knots of minimal warping degree two are the following 9 knots: 
$5_1, 5_2, 6_1, 6_2, 6_3, 7_6, 7_7, 8_{12}, 8_{18}$. 
\label{prime-alt-two}
\end{theorem}
\phantom{x}

\begin{figure}[ht]
\begin{center}
\includegraphics[width=70mm]{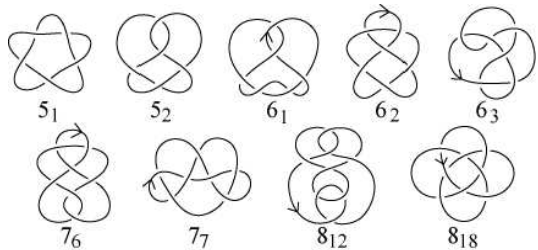}
\caption{Minimal diagrams of all the prime alternating knots with minimal warping degree two. The diagrams of $5_1$ and $5_2$ have warping degree two with both orientations. }
\label{mwd2kd}
\end{center}
\end{figure}

\noindent (See Figure \ref{mwd2kd}.) 
We remark that there are still minimal warping degree two knots for non-prime or non-alternating knots. 
For example, the Granny knot has minimal warping degree two (see Proposition \ref{16KPs}). 
We can also see that $md(8_{20})=md(8_{21})=2$ from Rolfsen's knot table. 

The rest of the paper is organized as follows. 
In Section 2, we define a half-curve and its length, which gives an estimation of the warping degree. 
In Section 3, we enumerate all the alternating knot diagrams of warping degree two and prove Theorem \ref{prime-alt-two}. 
In Section 4, we give a lower and upper bounds of the warping degree of some knot projections considering the ``r-factor''. 

\section{Half-curve and the warping degree}

\begin{figure}[ht]
\begin{center}
\includegraphics[width=25mm]{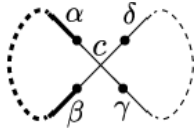}
\caption{Four points around a crossing $c$. Broken curves imply the connection. The half-curve $C^{\alpha}_{\beta}$ is thicken. }
\label{alpha}
\end{center}
\end{figure}

\begin{figure}[ht]
\begin{center}
\includegraphics[width=70mm]{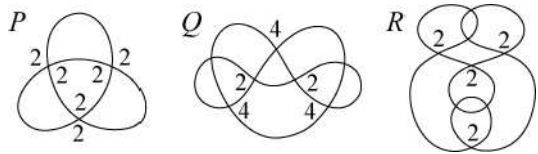}
\caption{Lengths of half-curves. The value of $l( C^{\alpha}_{\beta} )$ is given near the endpoints of each half-curve $C^{\alpha}_{\beta}$. The length is not defined at each empty corner. }
\label{l-ex}
\end{center}
\end{figure}

\noindent In this section, we define a half-curve and its length for knot projections and knot diagrams to give an estimation for the warping degree. 
For a knot projection $P$, give four points $\alpha, \beta, \gamma$ and $\delta$ around a crossing $c$ as shown in Figure \ref{alpha}. 
We call the curve from $\alpha$ to $\beta$ which does not include $c$ a {\it half-curve based on $c$}, and denote it by $C^{\alpha}_{\beta}$. 
If a half-curve $C^{\alpha}_{\beta}$ does not have a self-crossing, we define the {\it length}, $l( C^{\alpha}_{\beta} )$, of $C^{\alpha}_{\beta}$ to be the number of crossings on $C^{\alpha}_{\beta}$, that is, the number of crossings between $C^{\alpha}_{\beta}$ and $C^{\gamma}_{\delta}$. 
If $C^{\alpha}_{\beta}$ has a self-crossing, the length is not defined. 
See Figure \ref{l-ex} for an example. 
We note that the value of the length is a non-negative even integer. 

\noindent For knot diagrams, we define a half-curve and its length as well as knot projections. 
We have the following lemma: 

\phantom{x}
\begin{lemma}
If an oriented alternating knot diagram $D$ has a half-curve $C^{\alpha}_{\beta}$ whose length is $n$, then $d(D) \geq n/2$. 
\label{length-wd}
\end{lemma}
\phantom{x}

\begin{proof}
Let $c$ be a crossing of $D$ which $C^{\alpha}_{\beta}$ is based on. 
As mentioned in \cite{shimizu-wd}, for any oriented alternating knot diagram $D$, the warping degree $d(D)$ is obtained by taking a base point anywhere just before an over-crossing, or just after an under-crossing. 
Take such a base point $b$ around $c$ as shown in Figure \ref{beta}. 
\begin{figure}[ht]
\begin{center}
\includegraphics[width=90mm]{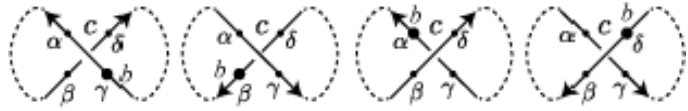}
\caption{Take a base point $b$ just before an over-crossing or just after an under-crossing of $c$. }
\label{beta}
\end{center}
\end{figure}
Then, a crossing $p$ on $C^{\alpha}_{\beta}$ is a warping crossing point of $D_b$ if $p$ is a crossing where $C^{\alpha}_{\beta}$ is under and $C^{\gamma}_{\delta}$ is over. 
Since $D$ is alternating, the half-curve $C^{\alpha}_{\beta}$ has the same number of over-crossings and under-crossings, and the number equals half the value of $l( C^{\alpha}_{\beta} )$. 
Recall that $C^{\alpha}_{\beta}$ has no self-crossings. 
Hence we have the inequality $d(D) \geq l( C^{\alpha}_{\beta} ) /2$. 
\end{proof}
\phantom{x}

For a knot projection $P$, we define the {\it length of $P$} to be the maximal value of the length for all half-curves of $P$ and denote it by $l(P)$. 
For example, the knot projections $P, Q$ and $R$ in Figure \ref{l-ex} have length two, four and two, respectively. 
A knot projection or knot diagram is {\it nontrivial} if it has at least one crossing. 
The following proposition implies that the length is well-defined for nontrivial knot projections. 

\phantom{x}
\begin{proposition}
Every nontrivial knot projection has a half-curve which has no self-crossings. 
\label{length-defined}
\end{proposition}

\phantom{x}
\begin{proof}
For a crossing $c$ of $P$, if a half-curve $C$ based on $c$ has a self-crossing $c'$, retake the half-curve $C'$ based on $c'$ which is included by $C$. 
If $C'$ has a self-crossing $c''$, similarly retake the half-curve $C''$ based on $c''$. 
Repeat this procedure and we obtain a half-curve without a self-crossing. 
\end{proof}
\phantom{x}

\noindent The length $l(D)$ of a knot diagram $D$ is well-defined as well as knot projections. 
From Lemma \ref{length-wd}, a lower bound for the warping degree of an oriented alternating knot diagram is given in the following corollary.

\phantom{x}
\begin{corollary}
For an oriented alternating nontrivial knot diagram $D$, the inequality $d(D) \geq l(D)/2$ holds. 
\label{dDlD}
\end{corollary}
\phantom{x}

\noindent We have the following corollary. 

\phantom{x}
\begin{corollary}
Any nontrivial alternating knot diagram $D$ with $d(D) \leq 2$ has length four, two or zero. 
\label{420}
\end{corollary}
\phantom{x}

\begin{proof}
The value of $l(D)$ should be four or less from $l(D)/2 \leq d(D) \leq 2$ by Corollary \ref{dDlD}, and be non-negative even integer by definition. 
\end{proof}
\phantom{x}

A {\it tangle} is a portion of a knot projection or knot diagram which can be bounded by a circle where arcs intersect transversely. 
We call each arc in a tangle which has the endpoints on the circle a {\it segment}. 
An {\it $n$-tangle} is a tangle of $n$ segments. 
We have the following. 

\phantom{x}
\begin{lemma}
Let $T$ be an $n$-tangle ($n \geq 2$) which is a portion of an alternating knot diagram. 
If two segments intersect mutually two or more times, then $T$ includes at least one warping crossing point by taking a base point anywhere outside $T$. 
\label{2tangle}
\end{lemma}
\phantom{x}

\begin{proof}
Let $S, R$ be two segments in $T$ which have two or more mutual crossings. 
Remark that $S$ and $R$ may also have self-crossings or crossings with other segments. 
Since the diagram is alternating, $S$ and $R$ have over-crossings and under-crossings alternatively. 
For each self-crossing, the segment has one over-crossing and one under-crossing. 
Hence, $S$ and $R$ have at least one non-self under-crossing because they have two or more non-self crossings. 
If we take a base point outside $T$, one of $S$ and $R$, which we encounter first from the base point, has a warping crossing point. 
\end{proof}
\phantom{x}

\noindent We have the following: 

\phantom{x}
\begin{corollary}
Let $T$ be an $n$-tangle which is a portion of an alternating knot diagram. 
If $T$ has a $k$-gon created by $k$ segments ($k \geq 2$), then $T$ includes a warping crossing point by taking a base point anywhere outside of $T$. 
\label{k-gonT}
\end{corollary}

\phantom{x}
\begin{proof}
Since the $k$ segments have a non-self under-crossings, one of them we encounter first has a warping crossing point when we take a base point outside $T$. 
\end{proof}
\phantom{x}

\noindent In this paper, we say that half-curves, tangles, polygons or regions are {\it disjoint} when they share no crossings. 
We have the following. 

\phantom{x}
\begin{corollary}
If an alternating knot diagram $D$ has disjoint $k_1$-gon, $k_2$-gon, $\dots$ and $k_n$-gon ($k_1 , k_2 , \dots , k_n \geq 2$) and a crossing $c$, then $d(D) \geq n$. 
\label{k1k2k3}
\end{corollary}
\phantom{x}

\begin{proof}
Assume that each $k_i$-gon is in a $k_i$-tangle $T_i$ which is obtained by taking a regular neighborhood of the $k_i$-gon. 
Take a base point $b$ just before an over-crossing of the crossing $c$. 
From Corollary \ref{k-gonT}, each $T_i$ has a warping crossing point, and therefore $D_b$ has at least $n$ warping crossing points. 
\end{proof}
\phantom{x}

\noindent We show the following corollary for knot projections which will be used in the next section. 

\phantom{x}
\begin{corollary}
Let $T$ be an $n$-tangle which is a portion of a knot projection $P$. 
If a segment $S$ of $T$ has a self-crossing, then $S$ includes at least one half-curve of $P$ which has no self-crossings. 
\label{tangle-hc}
\end{corollary}
\phantom{x}

\begin{proof}
Let $Q$ be a knot projection obtained from $S$ by connecting the endpoints outside of $T$. 
By Proposition \ref{length-defined}, $Q$ has a half-curve which has no self-crossings. 
Then, $S$ has the corresponding half-curve in $P$, and it has no self-crossings even if it has extra non-self-crossings. 
\end{proof}
\phantom{x}

\section{Construction of knot projections of warping degree two}

\begin{figure}[ht]
\begin{center}
\includegraphics[width=80mm]{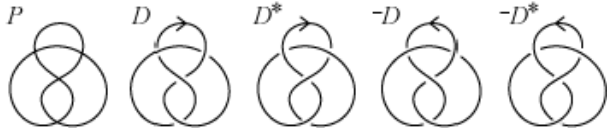}
\caption{From $d(D)=d(-D^*)=1$ and $d(-D)=d(D^*)=2$, the warping degree of the knot projection $P$ is 1. }
\label{41}
\end{center}
\end{figure}

The aim of this section is to list up all the reduced alternating knot diagrams $D$ with $d(D)=2$ and prove Theorem \ref{prime-alt-two}. 
The {\it warping degree $d(P)$ of a knot projection $P$} is the minimal value of warping degree for both alternating knot diagrams $D$ with $|D|=P$ and for both orientations. 
For example, the knot projection $P$ in Figure \ref{41} has $d(P)=1$. 
In \cite{shimizu-wd}, the equality $d(-D)=d(D^*)$ is shown, where $-D$ is the knot diagram $D$ with orientation reversed and $D^*$ is $D$ with all the crossings changed. 
Hence it is sufficient to check the warping degrees $d(D)$ and $d(D^*)$ for a fixed orientation to obtain the value of $d(P)$ for $P=|D|$. 
We have the following corollary from Corollary \ref{k1k2k3}.

\phantom{x}
\begin{corollary}
If a knot projection $P$ has disjoint $k_1$-gon, $k_2$-gon, $\dots$ and $k_n$-gon $(k_1, k_2, \dots ,k_n \geq 2)$ and a crossing, then $d(P) \geq n$. 
\label{k1k2k3p}
\end{corollary}
\phantom{x}

\noindent A knot projection $P$ is said to be {\it reducible} if $P$ has a region which meets itself diagonally around a crossing. 
A knot projection $P$ is said to be {\it reduced} if $P$ is not reducible. 
In this section, we consider reduced projections only, and then we do not need to consider half-curves with length zero. 
Let $P$ be a reduced alternating knot projection with $d(P)=2$. 
Then, we have the necessary condition that the length of $P$ is four or two by Corollary \ref{420}, and we have the following three cases by Lemma \ref{length-wd} and Corollary \ref{k1k2k3p}. 
Note that a half-curve of length two creates a 2-gon. \\

\noindent {\bf Case 1. } $l(P)=4$ and $P$ has exactly one half-curve of length four and does not have other disjoint half-curves. \\
{\bf Case 2. } $l(P)=2$ and $P$ has exactly two disjoint half-curves with length two, \\
{\bf Case 3. } $l(P)=2$ and $P$ does not have two or more disjoint half-curves with length two. \\

\begin{figure}[ht]
\begin{center}
\includegraphics[width=60mm]{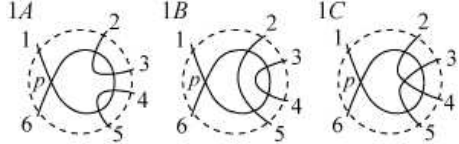}
\caption{Tangles $1A, 1B$ and $1C$ including a half-curve of length four. }
\label{case1}
\end{center}
\end{figure}

\noindent {\bf Case 1}: 
In this case, $P$ has one of the three tangles $1A, 1B$ and $1C$ in Figure \ref{case1} as a half-curve of length four in the following reason. 
The segment between $1$ and $6$ in each tangle includes a half-curve of length four. 
The other two segments may have at most one mutual crossing and have no self-crossings inside the half-curve by Lemma \ref{2tangle} and Corollary \ref{tangle-hc} by considering a 2-tangle inside the half-curve. 
Thus we have exactly three cases, $1A, 1B$ and $1C$. 
To list up knot projections of warping degree two, we connect the endpoints by three segments outside of each tangle satisfying the following necessary conditions for $d(P)=2$: \\

\noindent (i) The result is a knot projection (not a link projection). \\
(ii) Any pair of segments have at most one mutual crossing (by Lemma \ref{2tangle}). \\
(iii) Any segment does not have a self-crossing (by Corollary \ref{tangle-hc}). \\

\noindent For each case of $1A, 1B$ and $1C$, there are eight ways of connection satisfying the condition (i); 
The endpoint of 1 has four choices of endpoints to connect, except 6. 
Then the endpoint 6 has two choices to connect, the endpoints of the segment which is not connected to 1. 
Thus we have eight connections. \\

\begin{figure}[ht]
\begin{center}
\includegraphics[width=80mm]{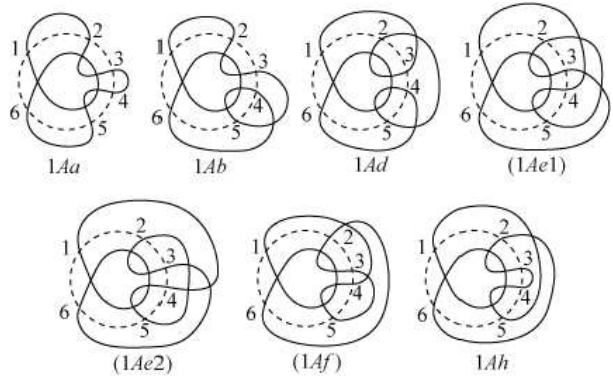}
\caption{All the knot projections $1Aa, 1Ab, 1Ad, 1Ae1, 1Ae2, 1Af$ and $1Ah$ satisfying the conditions (ii) and (iii) obtained by the connections $a, b, d, e, f$ and $h$ in Case $1A$. }
\label{case1A}
\end{center}
\end{figure}

For Case $1A$, the eight combinations satisfying the condition (i) are the following $a$ to $h$: \\
$a (1-2, \ 3-4, \ 5-6)$, $b (1-2, \ 3-5, \ 4-6)$, $c (1-3, \ 2-4, \ 5-6)$, $d (1-3, \ 2-5, \ 4-6)$, $e (1-4, \ 2-5, \ 3-6)$, $f (1-4, \ 2-6, \ 3-5)$, $g (1-5, \ 2-4, \ 3-6)$, $h (1-5, \ 2-6, \ 3-4)$ \\
Here, two combinations such that one is obtained from the other by subtracting all the numbers from 7 are the horizontal-reflection each other. 
For example, by subtracting from 7 for all the numbers of $b (1-2, \ 3-5, \ 4-6)$, we obtain $(6-5, \ 4-2, \ 3-1)$, and this combination is equivalent to $c$. 
The combinations $f$ and $g$ are also a horizontal-reflection pair. 
Then we can delete $c$ and $g$ from the list. 
In Figure \ref{case1A}, the knot projections which is a result of the combinations $a, b, d, e, f$ and $h$ satisfying the conditions (ii) and (iii), $1Aa$, $1Ab$, $1Ad$, $1Ae1$, $1Ae2$, $1Af$, $1Ah$, are shown. 
The knot projections $1Ae1$ and $1Ae2$ have a 3-gon out of the tangle, and have warping degree three or more when we take a base point at the crossing where the length-four half-curve is based on by Corollary \ref{k-gonT}. 
The knot projection $1Af$ has warping degree three, and the other knot projections have two. \\

\begin{figure}[ht]
\begin{center}
\includegraphics[width=80mm]{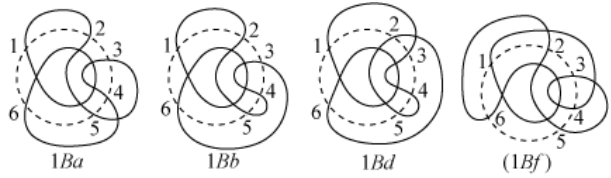}
\caption{All the knot projections $1Ba, 1Bb, 1Bd$ and $1Bf$ satisfying the conditions (ii) and (iii) obtained by the connections $a, b, d$ and $f$ in Case $1B$. }
\label{case1B}
\end{center}
\end{figure}

For Case $1B$, the eight combinations satisfying the condition (i) are the following $a$ to $h$: \\
$a (1-2, \ 3-5, \ 4-6)$, $b (1-2, \ 3-6, \ 4-5)$, $c (1-3, \ 2-4, \ 5-6)$, $d (1-3, \ 2-6, \ 4-5)$, $e (1-4, \ 2-3, \ 5-6)$, $f (1-4, \ 2-6, \ 3-5)$,  $g (1-5, \ 2-3, \ 4-6)$, $h (1-5, \ 2-4, \ 3-6)$ \\
Here, $a$ and $c$, $b$ and $e$, $d$ and $g$, and $f$ and $h$ are the horizontal-reflection pairs. 
In Figure \ref{case1B}, the knot projections are shown, and we can see that the knot projection $1Bf$ has warping degree three and the others have two. \\

\begin{figure}[ht]
\begin{center}
\includegraphics[width=80mm]{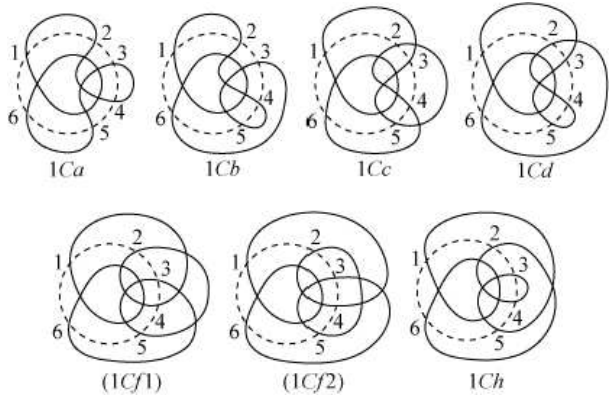}
\caption{All the knot projections $1Ca, 1Cb, 1Cc, 1Cd, 1Cf1, 1Cf2$ and $1Ch$ satisfying the conditions (ii) and (iii) obtained by the connections $a, b, c, d$ and $h$ in Case $1C$. }
\label{case1C}
\end{center}
\end{figure}

For Case $1C$, the eight combinations satisfying the condition (i) are the following $a$ to $h$: \\
$a (1-2, \ 3-4, \ 5-6)$, $b (1-2, \ 3-6, \ 4-5)$, $c (1-3, \ 2-5, \ 4-6)$, $d (1-3, \ 2-6, \ 4-5)$, $e  (1-4, \ 2-3, \ 5-6)$, $f (1-4, \ 2-5, \ 3-6)$, $g (1-5, \ 2-3, \ 4-6)$, $h (1-5, \ 2-6, \ 3-4)$\\
Here, $b$ and $e$, and $d$ and $g$ are the horizontal-reflection pairs. 
In Figure \ref{case1C}, the knot projections are shown, and we can see that $1Cf1$ and $1Cf2$ have a 3-gon out of the tangle, and the other knot projections have warping degree two. \\

\begin{figure}[ht]
\begin{center}
\includegraphics[width=40mm]{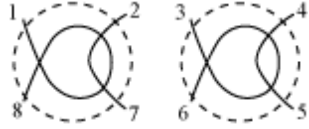}
\caption{Two tangles of half-curves of length two. }
\label{case2}
\end{center}
\end{figure}

\noindent {\bf Case 2}: 
In this case, $P$ has the two disjoint tangles in Figure \ref{case2}. 
The outer of the two discs bounded by broken circles forms an annulus, and to list up knot projections, we connect the endpoints 1 to 8 by four segments on the annulus satisfying: \\

\noindent (i) The result is a reduced knot projection. \\
(ii) Any segment does not have a self-crossing. \\
(iii) Any segment has at most two crossings on it. \\
(iv) Any pair of segments have at most one mutual crossing. \\

\begin{figure}[ht]
\begin{center}
\includegraphics[width=120mm]{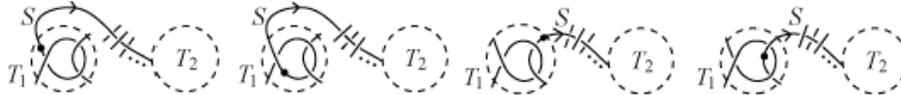}
\caption{If a segment $S$ has three or more non-self crossings, the warping degree is greater than two. }
\label{Sthree}
\end{center}
\end{figure}

\noindent We have the condition (ii) because if a segment $S$ has a self-crossing, then the annulus includes a half-curve with a length similarly to Corollary \ref{tangle-hc} and contradicts the situation of Case 2. 
We also have the condition (iii) in the following reason. 
Consider an oriented alternating knot diagram $D$ with $|D|=P$. 
If a segment $S$ on the annulus has three or more non-self crossings, take a base point in one of the 2-tangles, say $T_1$, just before an over-crossing or just after an under-crossing so that the first segment from the base point is $S$. 
There are essentially four types shown in Figure \ref{Sthree}. 
In the three cases from the left, the two 2-tangles have warping crossing points because they have a 2-gon (Corollary \ref{k-gonT}), and $S$ has one or more warping crossing points. 
In the other case, the 2-tangle $T_2$ which is not $T_1$ has a warping crossing point, and $S$ has two or more warping crossing points. 
Hence, if $S$ has three or more crossings, then the warping degree is three or more. 
We have the condition (iv) from Lemma \ref{2tangle} and Corollary \ref{k1k2k3p}.

We have 48 combinations satisfying the condition (i) because the endpoint of 1 has six choices to connect, except the end point of 8. 
Then the endpoint of 8 has four choices to connect, the endpoints of the rest two segments. 
And then the rest four endpoints have two combinations. 
From a combination, we obtain the horizontal-reflected combination by subtracting all the numbers from 9. 
We have the following 24 combinations $A$ to $X$ up to horizontal-reflections. \\

\begin{align*}
A: & \ (1-2, \ 3-4, \ 5-7, \ 6-8), \ (1-3, \ 2-4, \ 5-6, \ 7-8) \\
B: & \ (1-2, \ 3-4, \ 5-8, \ 6-7), \ (1-4, \ 2-3, \ 5-6, \ 7-8) \\
C: & \ (1-2, \ 3-5, \ 4-7, \ 6-8), \ (1-3, \ 2-5, \ 4-6, \ 7-8) \\
D: & \ (1-2, \ 3-5, \ 4-8, \ 6-7), \ (1-5, \ 2-3, \ 4-6, \ 7-8) \\
E: & \ (1-2, \ 3-7, \ 4-6, \ 5-8), \ (1-4, \ 2-6, \ 3-5, \ 7-8) \\
F: & \ (1-2, \ 3-7, \ 4-8, \ 5-6), \ (1-5, \ 2-6, \ 3-4, \ 7-8) \\
G: & \ (1-2, \ 3-8, \ 4-6, \ 5-7), \ (1-6, \ 2-4, \ 3-5, \ 7-8) \\
H: & \ (1-2, \ 3-8, \ 4-7, \ 5-6), \ (1-6, \ 2-5, \ 3-4, \ 7-8) \\
I: & \ (1-3, \ 2-4, \ 5-8, \ 6-7), \ (1-4, \ 2-3, \ 5-7, \ 6-8) \\
J: & \ (1-3, \ 2-5, \ 4-8, \ 6-7), \ (1-5, \ 2-3, \ 4-7, \ 6-8) \\
K: & \ (1-3, \ 2-6, \ 4-7, \ 5-8), \ (1-4, \ 2-5, \ 3-7, \ 6-8) \\
L: & \ (1-3, \ 2-6, \ 4-8, \ 5-7), \ (1-5, \ 2-4, \ 3-7, \ 6-8) \\
M: & \ (1-3, \ 2-8, \ 4-6, \ 5-7), \ (1-7, \ 2-4, \ 3-5, \ 6-8) \\
N: & \ (1-3, \ 2-8, \ 4-7, \ 5-6), \ (1-7, \ 2-5, \ 3-4, \ 6-8) \\
O: & \ (1-4, \ 2-5, \ 3-8, \ 6-7), \ (1-6, \ 2-3, \ 4-7, \ 5-8) \\
P: & \ (1-4, \ 2-6, \ 3-8, \ 5-7), \ (1-6, \ 2-4, \ 3-7, \ 5-8) \\
Q: & \ (1-4, \ 2-8, \ 3-5, \ 6-7), \ (1-7, \ 2-3, \ 4-6, \ 5-8) \\
R: & \ (1-4, \ 2-8, \ 3-7, \ 5-6), \ (1-7, \ 2-6, \ 3-4, \ 5-8) \\
S: & \ (1-5, \ 2-4, \ 3-8, \ 6-7), \ (1-6, \ 2-3, \ 4-8, \ 5-7) \\
T: & \ (1-5, \ 2-6, \ 3-8, \ 4-7), \ (1-6, \ 2-5, \ 3-7, \ 4-8) \\
U: & \ (1-5, \ 2-8, \ 3-4, \ 6-7), \ (1-7, \ 2-3, \ 4-8, \ 5-6) \\
V: & \ (1-5, \ 2-8, \ 3-7, \ 4-6), \ (1-7, \ 2-6, \ 3-5, \ 4-8) \\
W: & \ (1-6, \ 2-8, \ 3-4, \ 5-7), \ (1-7, \ 2-4, \ 3-8, \ 5-6) \\
X: & \ (1-6, \ 2-8, \ 3-5, \ 4-7), \ (1-7, \ 2-5, \ 3-8, \ 4-6)
\end{align*}

\begin{figure}[ht]
\begin{center}
\includegraphics[width=80mm]{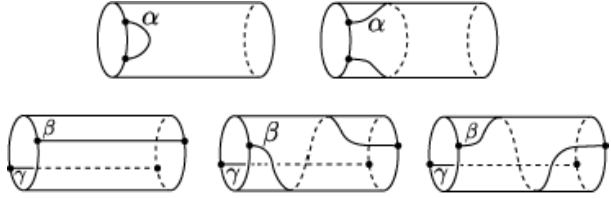}
\caption{A segment $\alpha$ whose endpoints are on the same side has two choices to place on the annulus. 
Two segments $\beta$ and $\gamma$ whose endpoints are on the both sides have relatively three choices. 
We say that $\beta$ is parallel to $\gamma$ for the lower left case, right-handed to $\gamma$ for the lower middle, and left-handed to $\gamma$ for the lower right. }
\label{annulus}
\end{center}
\end{figure}

\noindent Next we connect the endpoints by four segments on the annulus for each combination considering the conditions (ii) to (iv). 
We note that a segment whose endpoints are on the same side has two choices to place on an annulus, and two segments whose endpoints are on the both sides have relatively three choices to place on annulus as shown in Figure \ref{annulus}. 
Recall that two segments can not intersect two or more times by condition (iv). \\

\begin{figure}[ht]
\begin{center}
\includegraphics[width=85mm]{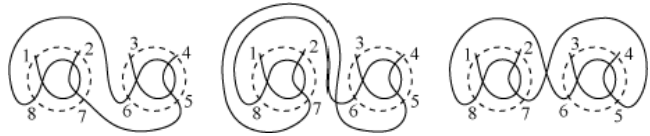}
\caption{The segment 5-7 with parallel, right-handed and left-handed to the segment 6-8.  }
\label{2A}
\end{center}
\end{figure}

\begin{figure}[ht]
\begin{center}
\includegraphics[width=120mm]{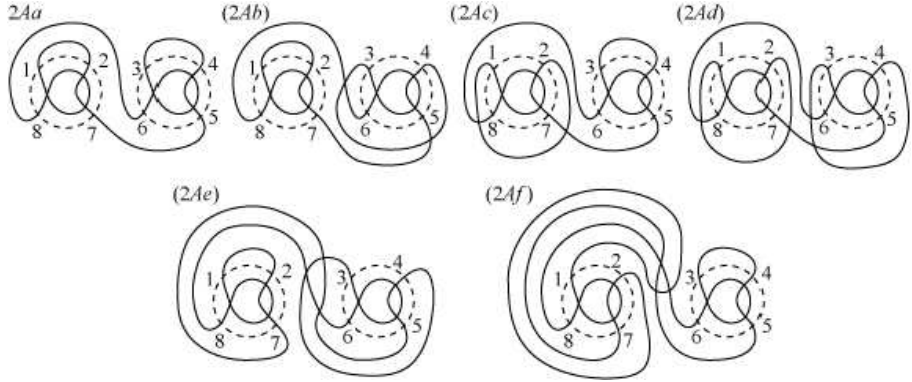}
\caption{All the knot projections satisfying the conditions (ii) to (iv) for case $2A$. Only the knot projection $2Aa$ has warping degree two. }
\label{case2A}
\end{center}
\end{figure}

For Case $2A (1-2, \ 3-4, \ 5-7, \ 6-8)$, fix the segment 6-8 without loss of generality. 
Then the segment 5-7 has three choices as shown in Figure \ref{2A}, parallel, right-handed and left-handed to 6-8. 
For the case that 5-7 is parallel to 6-8, the segment 1-2 has two choices with the condition (iv) that 1-2 and any segment can not intersect twice. 
In the same way, 3-4 has two choices. 
Thus we obtain four knot projections $2Aa, 2Ab, 2Ac$ and $2Ad$ shown in the upper in Figure \ref{case2A}. 
The knot projections $2Ab, 2Ac, 2Ad$ have three disjoint 2-gons and a crossing, and therefore have warping degree three or more by Corollary \ref{k1k2k3p}. 
The knot projection $2Aa$ has warping degree two. 

Next, for the case that 5-7 is right-handed to 6-8, we need to place 1-2 and 3-4 so that the result is a reduced knot projection. 
We obtain two knot projections satisfying the conditions (ii) to (iv), $2Ae$ and $2Af$, shown in Figure \ref{case2A}. 
They have a 3-gon on the annulus, and have warping degree three or more. 
For the case that 5-7 is left-handed to 6-8, there are no way to make it reduced with the conditions (ii) to (iv). \\

\begin{figure}[ht]
\begin{center}
\includegraphics[width=85mm]{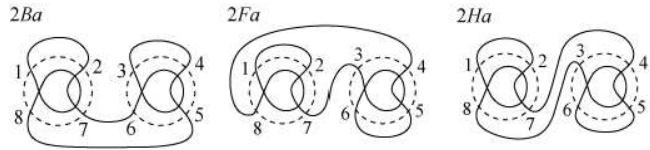}
\caption{All the knot projections with warping degree two for Cases $2B, 2F$ and $2H$. }
\label{case2B}
\end{center}
\end{figure}

In the same way to Case $2A$, we obtain knot projections $2Ba, 2Fa$ and $2Ha$, shown in figure \ref{case2B}, satisfying the conditions (ii) - (iv) for Cases $2B, 2F$ and $2H$, respectively. 
We can see that the warping degrees of $2Ba, 2Fa$ and $2Ha$ are all two. \\

\begin{figure}[ht]
\begin{center}
\includegraphics[width=120mm]{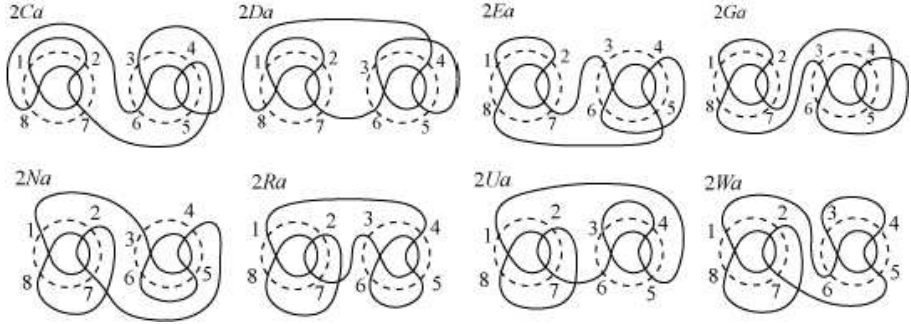}
\caption{All the knot projections with warping degree two for Cases $2C, 2D, 2E, 2G, 2N, 2R, 2U$ and $2W$. }
\label{case2C}
\end{center}
\end{figure}

For Case $2C (1-2, \ 3-5, \ 4-7, \ 6-8)$, there are three choices to place 4-7 and 6-8, similarly to the situation of Figure \ref{2A}. 
For the two cases that 4-7 is right-handed or left-handed to 6-8, i.e., 4-7 and 6-8 have a crossing, it can not avoid that the result is reducible or has a 3-gon on the annulus with the conditions (ii) to (iv), similarly to Case $2A$. 
Hence 4-7 and 6-8 should be parallel. 
The segments 1-2 and 3-5 have two choices, and one of them creates a 2-gon which is disjoint to the 2-gons in the tangles. 
Hence 1-2 and 3-5 should be on the side which does not create a 2-gon. 
Thus we obtain a knot projection $2Ca$, shown in Figure \ref{case2C}, which has warping degree two. \\

In the same way to Case $2C$, we obtain one knot projection, shown in Figure \ref{case2C}, satisfying the conditions and without three disjoint 2-gons for Cases $2D$, $2E$, $2G$, $2N$, $2R$, $2U$ and $2W$. 
All the knot projections have warping degree two. \\

\begin{figure}[ht]
\begin{center}
\includegraphics[width=25mm]{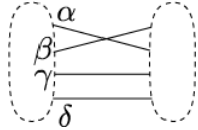}
\caption{The segments $\alpha , \beta , \gamma$ and $\delta$ for Cases $2I, 2J, 2K, 2P, 2S$ and $2T$. }
\label{abcd}
\end{center}
\end{figure}

\begin{figure}[ht]
\begin{center}
\includegraphics[width=70mm]{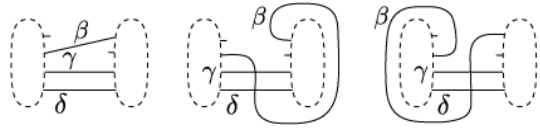}
\caption{For the case that $\gamma$ is parallel to $\delta$, $\beta$ has three choices, parallel, right-handed and left-handed to $\delta$. }
\label{pprl}
\end{center}
\end{figure}

For Case $2I (1-3, \ 2-4, \ 5-8, \ 6-7)$, we can see that every sugment has endpoints on both sides. 
In particular, the four segments can be placed on the annulus with exactly one crossing. 
We call segments $\alpha , \beta , \gamma , \delta$, where $\alpha$ and $\beta$ have a crossing at the minimal placement as shown in Figure \ref{abcd}. 
Now we enumerate all the possible placements of the segments with the conditions (ii) to (iv). 
Fix $\delta$ without loss of generality. 
Then $\gamma$ has three choices, parallel, right-handed or left-handed to $\delta$. 
For the case that $\gamma$ is parallel to $\delta$, then $\beta$ has three choices as well (see Figure \ref{pprl}). 
For the case that $\beta$ is parallel to $\delta$, then $\alpha$ has three choices, too. 
We can place $\alpha$ with the conditions for the case that $\alpha$ is parallel or left-handed to $\delta$, whereas we can not for the right-handed case because $\alpha$ and $\beta$ cross twice. 
For the case that $\beta$ is right-handed to $\delta$, we can place $\alpha$ only with parallel to $\delta$, otherwise $\alpha$ has three crossings. 
For the case that $\beta$ is left-handed to $\delta$, $\alpha$ and $\beta$ have two mutual crossings or some segments have three or more crossings for the three choices of $\alpha$. 

\begin{figure}[ht]
\begin{center}
\includegraphics[width=50mm]{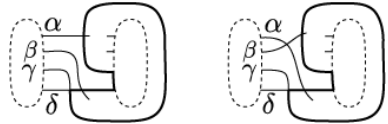}
\caption{When $\gamma$ is right-handed to $\delta$, $\alpha$ and $\beta$ must intersect with $\gamma$ or $\delta$. 
Since each segment can have at most two crossings, we have the two situations, $\alpha$ intersects $\gamma$ and $\beta$ intersects $\delta$, and $\alpha$ intersects $\delta$ and $\beta$ intersects $\gamma$. }
\label{r-ab}
\end{center}
\end{figure}

For the case that $\gamma$ is right-handed to $\delta$, the segments $\alpha$ and $\beta$ must intersect $\gamma$ or $\delta$ since the left-side endpoints of $\alpha$ and $\beta$ are bounded by $\gamma$, $\delta$ and the boundary of the right-side tangle. 
Since $\gamma$ and $\delta$ already have a crossing, they can have at most one crossing with $\alpha$ or $\beta$. 
Hence we have two cases as shown in Figure \ref{r-ab}. 
We can place $\alpha$ and $\beta$ satisfying the conditions for the left-hand case, and we can not for the right-hand case without crossing twice. 

\begin{figure}[ht]
\begin{center}
\includegraphics[width=110mm]{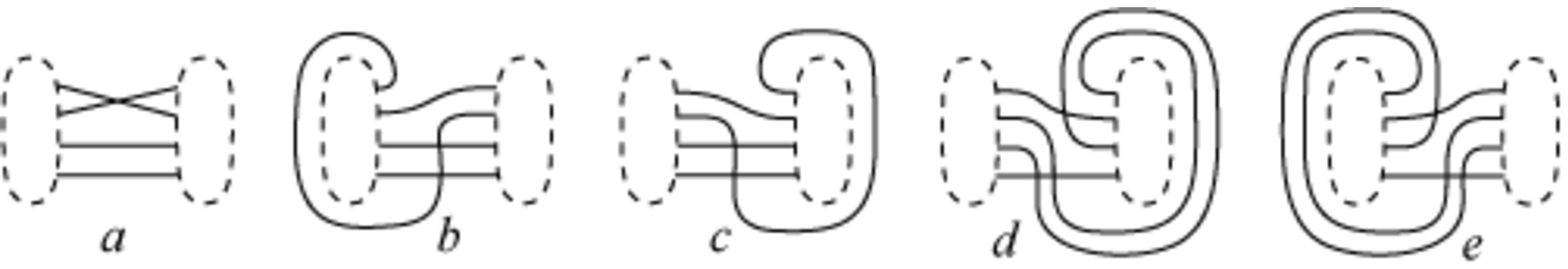}
\caption{The five placements $a$ to $e$ of the segments $\alpha , \beta , \gamma$ and $\delta$ satisfying the conditions (ii) to (iv). }
\label{abcde}
\end{center}
\end{figure}

\begin{figure}[ht]
\begin{center}
\includegraphics[width=125mm]{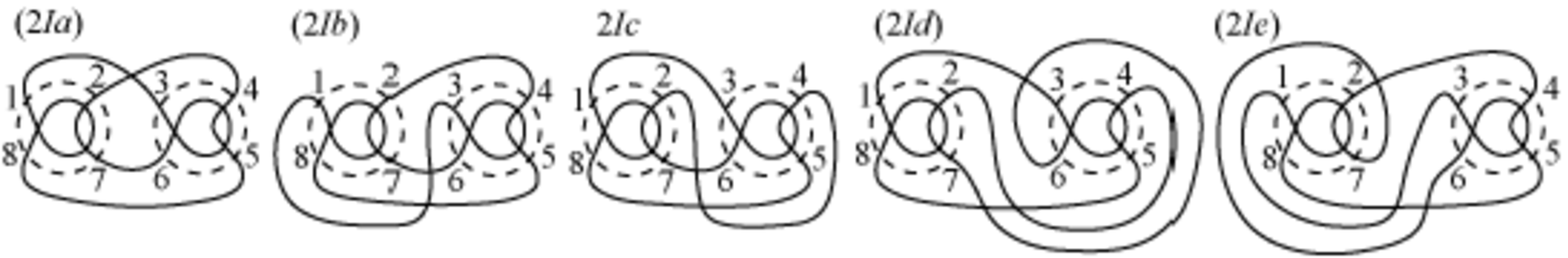}
\caption{The knot projections obtained by the placements $a$ to $e$ for Case $2I$. 
The knot projection $2Ia$ belongs to Case 1. 
Only the knot projection $2Ic$ has warping degree two. }
\label{case2I}
\end{center}
\end{figure}

Finally, for the case that $\gamma$ is left-handed to $\delta$, we obtain one knot projection in the same way to the case that $\gamma$ is right-handed to $\delta$. 
Thus, we have the five placements of $\alpha , \beta , \gamma$ and $\delta$, shown in Figure \ref{abcde}, and therefore we have five knot projections $2Ia, 2Ib, 2Ic, 2Id$ and $2Ie$, shown in Figure \ref{case2I} for case $2I$. 
The knot projection $2Ia$ has a half-curve of length four, and $2Ib$ and $2Id$ have three disjoint 2-gons. 
The knot projection $2Ie$ has warping degree four, and $2Ic$ has two. \\

\begin{figure}[ht]
\begin{center}
\includegraphics[width=125mm]{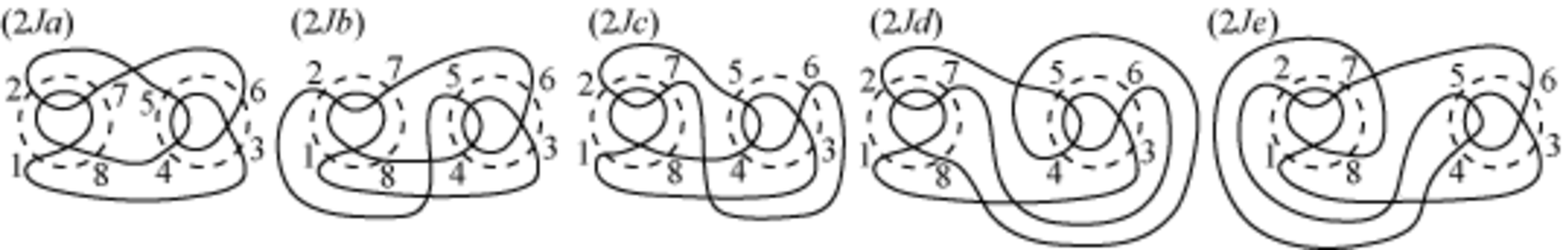}
\caption{All the knot projections obtained by the placements $a$ to $e$ for Case $2J$. 
The knot projection $2Ja$ belongs to Case 1. 
No knot projections have warping degree two except $2Ja$. }
\label{case2J}
\end{center}
\end{figure}

For Case $2J (1-3, \ 2-5, \ 4-8, \ 6-7)$, we obtain the five knot projections, shown in Figure \ref{case2J} from the five connections shown in Figure \ref{abcde}. 
The knot projection $2Ja$ has a half-curve of length four. 
The knot projections $2Jc$ and $2Jd$ have three disjoint 2-gons. 
The knot projection $2Jb$ and $2Je$ have warping degree three. \\

\begin{figure}[ht]
\begin{center}
\includegraphics[width=125mm]{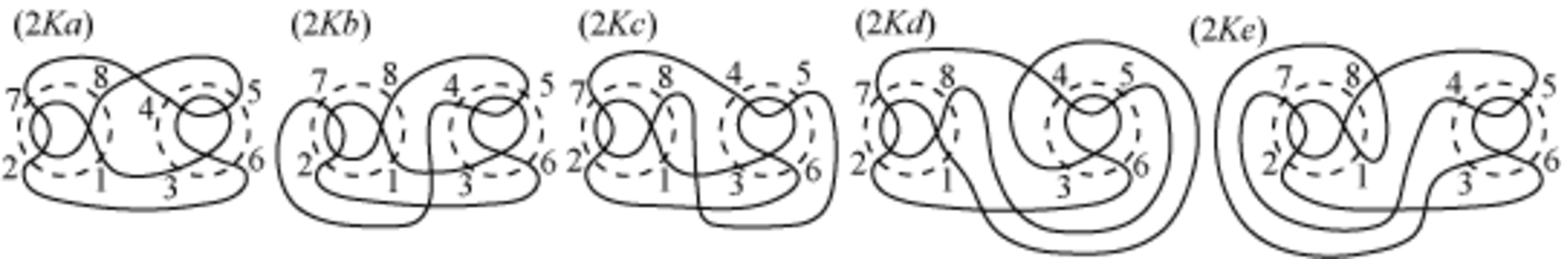}
\caption{All the knot projections obtained by the placements $a$ to $e$ for Case $2K$. 
The knot projection $2Ka$ belongs to Case 1. 
No knot projections have warping degree two except $2Ka$.  }
\label{case2K}
\end{center}
\end{figure}

For Case $2K (1-3, \ 2-6, \ 4-7, \ 5-8)$, we obtain the five knot projections, shown in Figure \ref{case2K} from the connections shown in Figure \ref{abcde}. 
The knot projection $2Ka$ has a half-curve of length four. 
The knot projections $2Kc$ and $2Ke$ have three disjoint 2-gons. 
The knot projections $2Kb$ and $2Kd$ have warping degree three. \\

\begin{figure}[ht]
\begin{center}
\includegraphics[width=125mm]{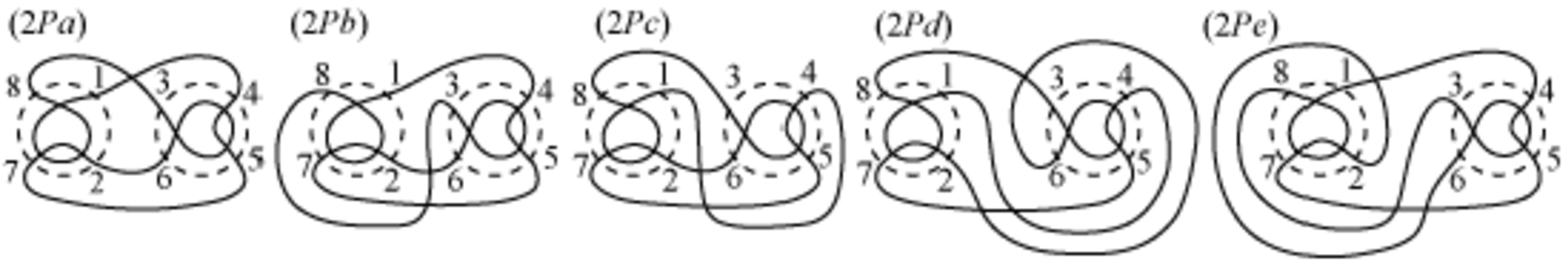}
\caption{All the knot projections obtained by the placements $a$ to $e$ for Case $2P$. 
No knot projections have warping degree two.  }
\label{case2P}
\end{center}
\end{figure}

For Case $2P (1-4, \ 2-6, \ 3-8, \ 5-7)$, we obtain the five knot projections, shown in Figure \ref{case2P} from the connections shown in Figure \ref{abcde}. 
The knot projections $2Pa, 2Pb$ and $2Pd$ have three disjoint 2-gons. 
The knot projection $2Pc$ has warping degree three, and $2Pe$ has four. \\

\begin{figure}[ht]
\begin{center}
\includegraphics[width=125mm]{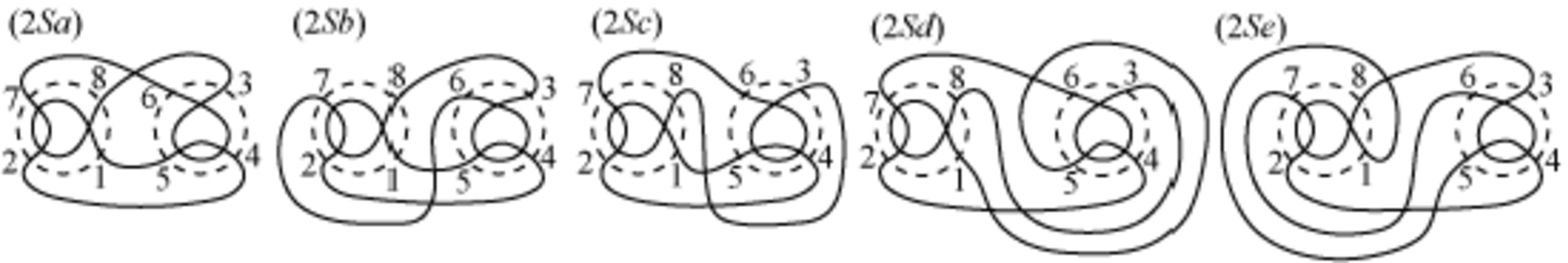}
\caption{All the knot projections obtained by the placements $a$ to $e$ for Case $2S$. 
No knot projections have warping degree two.  }
\label{case2S}
\end{center}
\end{figure}

For Case $2S (1-5, \ 2-4, \ 3-8, \ 6-7)$, we obtain the five knot projections, shown in Figure \ref{case2S} from the connections shown in Figure \ref{abcde}. 
The knot projections $2Sa, 2Sc$ and $2Se$ have three disjoint 2-gons. 
The knot projection $2Sb$ has warping degree three, and $2Sd$ has four. \\

\begin{figure}[ht]
\begin{center}
\includegraphics[width=125mm]{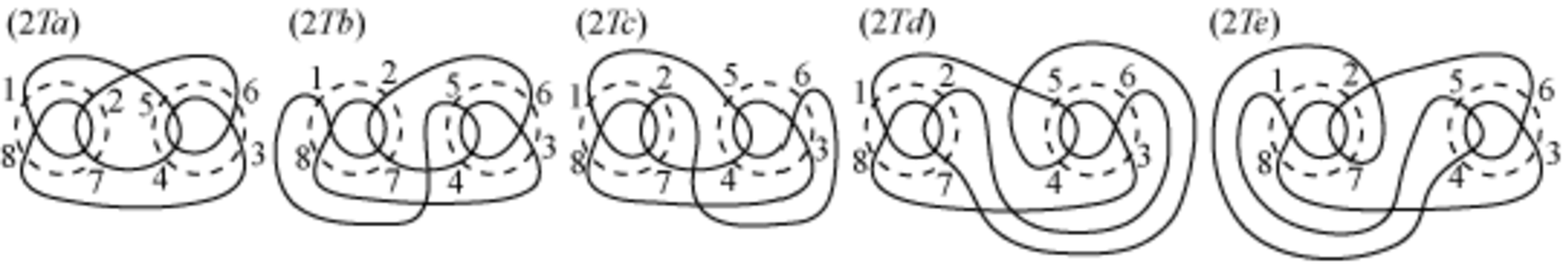}
\caption{All the knot projections obtained by the placements $a$ to $e$ for Case $2T$. 
The knot projection $2Ta$ belongs to Case 1. 
No knot projections have warping degree two except $2Ta$.  }
\label{case2T}
\end{center}
\end{figure}

For Case $2T (1-5, \ 2-6, \ 3-8, \ 4-7)$, we obtain the five knot projections, shown in Figure \ref{case2T} from the connections shown in Figure \ref{abcde}. 
The knot projection $2Ta$ has a half-curve of length four, and the other knot projections have three disjoint 2-gons. \\

\begin{figure}[ht]
\begin{center}
\includegraphics[width=120mm]{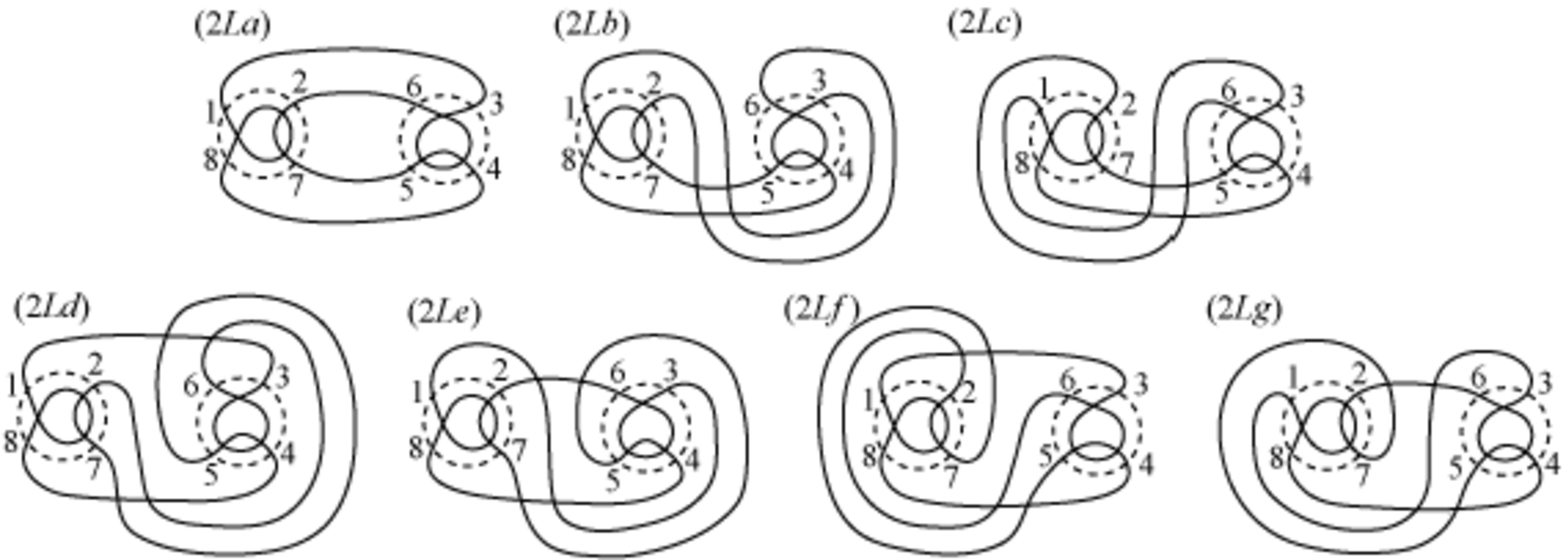}
\caption{All the knot projections with possible placements for Case $2L$. 
The knot projection $2La$ belongs to Case 1. 
No knot projections have warping degree two except $2La$.  }
\label{case2L}
\end{center}
\end{figure}

For Case $2L (1-3, \ 2-6, \ 4-8, \ 5-7)$, we obtain the seven knot projections, shown in Figure \ref{case2L}, in the same way to Case $2I$. 
More throughly, let $\delta , \gamma , \beta , \alpha$ be the segments 4-8, 5-7, 2-6, 1-3, respectively. 
Fix $\delta$ without loss of generality. 
We denote by $\delta (\gamma )=p, r, l$ when a segment $\gamma$ is parallel, right-handed, left-handed to $\delta$, respectively. 
We represent the status of $\gamma , \beta$ and $\alpha$ to $\delta$ by $( \delta (\gamma ), \delta (\beta ), \delta (\alpha ))$. 
When $( \delta (\gamma ), \delta (\beta ), \delta (\alpha )) = (p,p,p), (p,r,r), (p,l,l), (r,r,p), (r,p,r), (l,l,p)$ or $(l,p,l)$, we can place the segments satisfying the conditions (ii) to (iv), as shown in Figure \ref{case2L}. 
For the other 20 statuses, it is inevitable that some segments have three or more crossings. 

The knot projection $2La$ has a half-curve of length four. 
The knot projections $2Lc, 2Ld$ and $2Lg$ have three disjoint 2-gons, and $2Lb$ and $2Lf$ have a 4-gon on the annulus. 
The knot projection $2Le$ has warping degree three. \\

\begin{figure}[ht]
\begin{center}
\includegraphics[width=110mm]{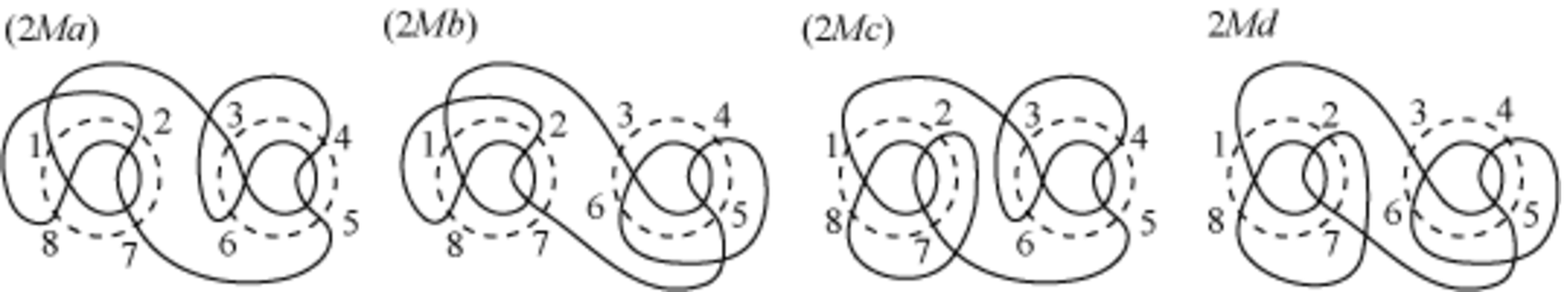}
\caption{All the knot projections satisfying the conditions for Case $2J$. 
Only the knot projection $2Md$ has warping degree two. }
\label{case2M}
\end{center}
\end{figure}

For Case $2M (1-3, \ 2-8, \ 4-6, \ 5-7)$, fix 5-7 without loss of generality. 
Then 1-3 has three choices. 
For the two cases that 1-3 and 5-7 have a crossing, it is inevitable that the result is reducible with the conditions. 
Hence 1-3 should be parallel to 5-7. 
Since 2-8 and 4-6 have two choices, we have four knot projections, $2Ma, 2Mb, 2Mc$ and $2Md$, shown in Figure \ref{case2M}. 
The knot projections $2Md$ has warping degree two. 
The others have three disjoint 2-gons. \\

\begin{figure}[ht]
\begin{center}
\includegraphics[width=120mm]{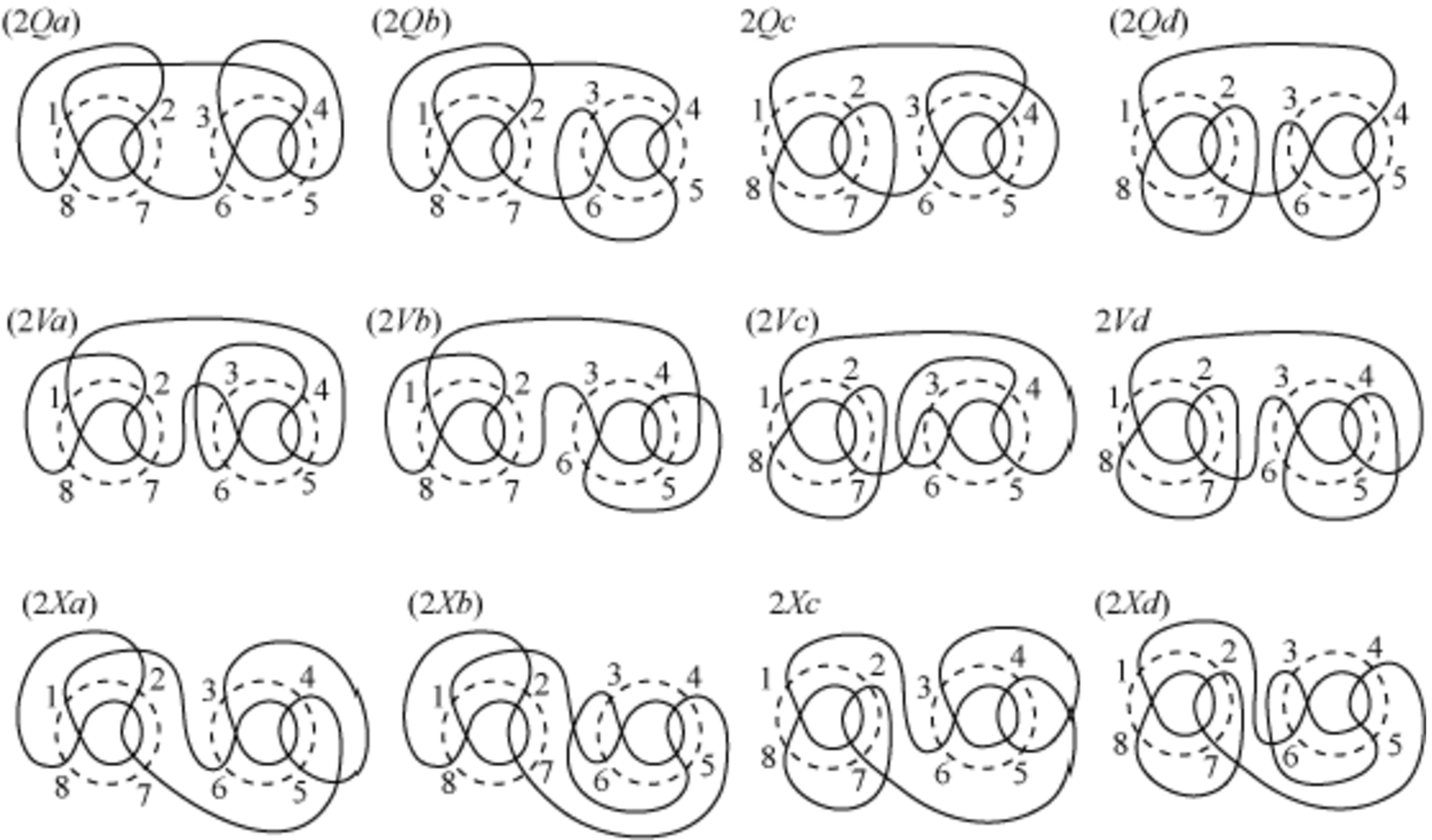}
\caption{All the knot projections satisfying the conditions for Cases $2Q, 2V$ and $2X$. 
Only the knot projections $2Qc, 2Vd$ and $2Xc$ have warping degree two. }
\label{case2Q}
\end{center}
\end{figure}

In the same way to Case $2M$, we obtain four knot projections, shown in Figure \ref{case2Q} for Cases $2Q, 2V$ and $2X$. 
The knot projections $2Qc, 2Vd, 2Xc$ have warping degree two, and the others have three disjoint 2-gons. \\

\begin{figure}[ht]
\begin{center}
\includegraphics[width=125mm]{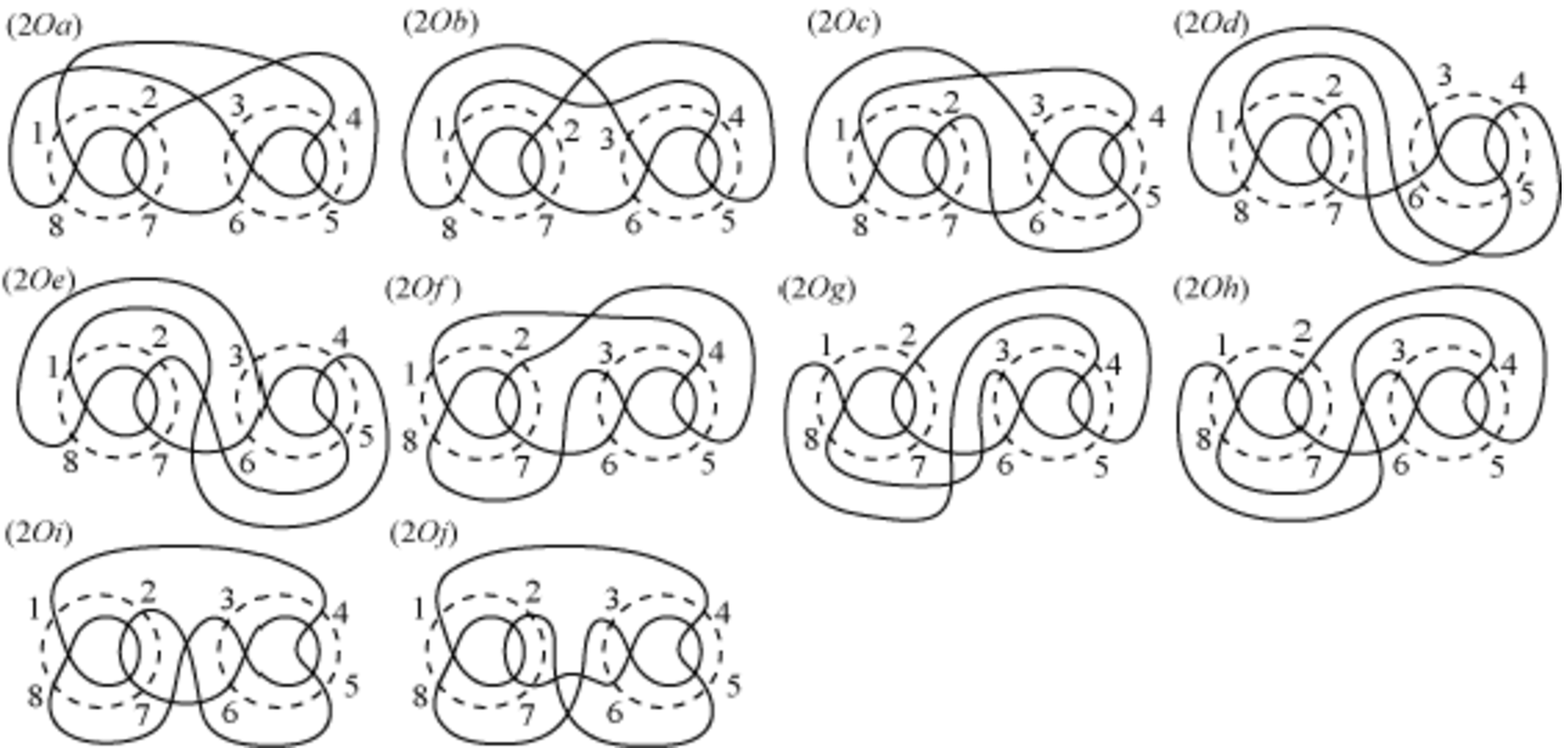}
\caption{All the knot projections satisfying the conditions for case $2O$. No knot projections have warping degree two. }
\label{case2O}
\end{center}
\end{figure}

For Case $2O (1-4, \ 2-5, \ 3-8, \ 6-7)$, fix 6-7 without loss of generality. 
Then 3-8 has three choices. 
For the case that 3-8 is parallel to 6-7, 2-5 has three choices. 
If 2-5 is left-handed to 6-7, then 2-5 has three or more crossings. 
Hence 2-5 should be parallel or right-handed to 6-7. 
Thus we obtain the knot projections $2Oa, 2Ob, 2Oc, 2Od$ and $2Oe$ shown in Figure \ref{case2O} for the case that 3-8 is parallel to 6-7. 
The knot projections $2Oa, 2Ob, 2Od$ and $2Oe$ have a 3-gon on the annulus. 
The knot projection $2Oc$ has warping degree three. 

For the case that 3-8 is right-handed to 6-7, if 2-5 is parallel to 6-7, then 3-8 has three or more crossings. 
If 2-5 is right-handed to 6-7, then 3-8 has three or more crossings with any placement of 1-4. 
If 2-5 if left-handed to 6-7, then 2-5 has three or more crossings. 

For the case that 3-8 is left-handed to 6-7, if 2-5 is left-handed to 6-7, then 2-5 have three or more crossings with any placement of 1-4. 
Hence 2-5 should be parallel or right-handed to 6-7. 
Thus we obtain the knot projections $2Of, 2Og, 2Oh, 2Oi, 2Oj$. 
The knot projection $2Of$ has three disjoint 2-gons, and the others have a 3-gon in the annulus. \\

\begin{figure}[ht]
\begin{center}
\includegraphics[width=13mm]{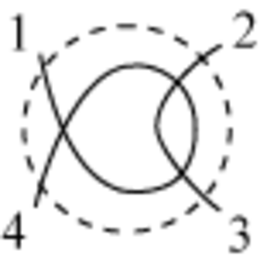}
\caption{For Case 3, a knot projection $P$ has a tangle of a half-curve of length two, and there are no disjoint tangles of a half-curve. }
\label{case3}
\end{center}
\end{figure}

\noindent {\bf Case 3}: 
In this case, $P$ has a tangle $T$ shown in Figure \ref{case3} which includes a half-curve of length two. 
The outer of the tangle $T$ forms a 2-tangle, say $\overline{T}$, satisfying the following conditions: \\

\noindent (i) The result of connecting is a reduced knot projection. \\
(ii) Any segment in $\overline{T}$ does not have a self-crossing. \\
(iii) Each segment in $\overline{T}$ has at most three crossings. \\

\noindent We have the condition (ii) to exclude Cases 1 and 2. 
Note that if a segment of $\overline{T}$ has a self-crossing, $\overline{T}$ includes a half-curve whose length is defined by Corollary \ref{tangle-hc}. 
We have condition (iii) in the same way to the condition (iii) in Case 2. 

\begin{figure}[ht]
\begin{center}
\includegraphics[width=90mm]{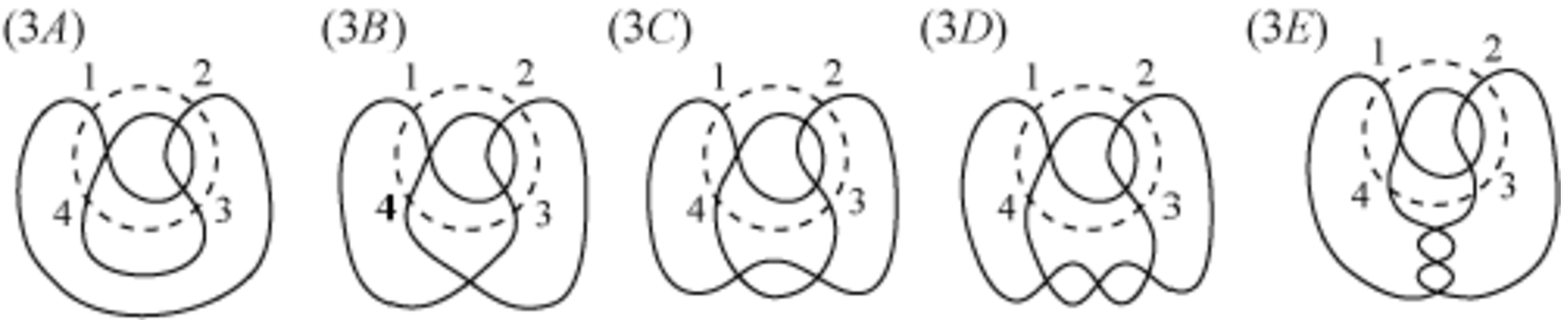}
\caption{All the knot projections obtained from the tangle in Figure \ref{case3} by connecting the segments satisfying the conditions for Case 3. 
The knot projections $3A$ and $3B$ have warping degree one. 
The other knot projections belong to Case 1. }
\label{3A}
\end{center}
\end{figure}

From the condition (i), we have the two connections ($1-2, \ 3-4$) and ($1-3, \ 2-4$). 
We have the four knot projections in Figure \ref{3A} which satisfy the conditions. 
The knot projections $3A$ and $3B$ have warping degree one. 
The knot projections $3C, 3D$ and $3E$ have a half-curve of length four. \\

\noindent In conclusion, we have the following. 

\phantom{x}
\begin{proposition}
All the reduced knot projections of warping degree two are the 16 knot projections shown in Figure \ref{mwd2kp}. 
\label{16KPs}
\end{proposition}
\phantom{x}

\noindent It is well-known that each reduced alternating diagram is a minimal diagram of an alternating knot, and furthermore each minimal diagram of a prime alternating knot is an alternating diagram as Tait's first conjecture, solved by Kauffman (\cite{kauffman}), Murasugi ({\cite{murasugi}) and Thistlethwaite (\cite{thistlethwaite}). 
We prove Theorem \ref{prime-alt-two}. 

\phantom{x}
\noindent {\it Proof of Theorem \ref{prime-alt-two}.} \ 
All the reduced knot projections of warping degree two are listed in Figure \ref{mwd2kp}. 
From the knot projection $1Aa$, we obtain a minimal diagram of $5_1$ as an alternating diagram $D$ with $|D|=1Aa$. 
From the knot projection $1Bb$, we obtain a minimal diagram of $5_2$. 
From the knot projection $1Ah$, we obtain a minimal diagram of $6_1$. 
From the knot projections $1Ab, 1Bd$ and $1Ca$, we obtain a minimal diagram of $6_2$. 
From the knot projections $1Ba$ and $1Cb$, we obtain a minimal diagram of $6_3$. 
From the knot projection $1Ch$, we obtain a minimal diagram of $7_6$. 
From the knot projections $1Ad$ and $1Cd$, we obtain a minimal diagram of $7_7$. 
From the knot projection $2Ic$, we obtain a minimal diagram of $8_{12}$. 
From the knot projection $1Cc$, we obtain a minimal diagram of $8_{18}$. 
From the other knot projections, we obtain minimal diagrams of composite knots. 
Hence, all the prime alternating knots of minimal warping degree two are the above nine knots. 
\hfill$\Box$

\phantom{x}

\begin{figure}[ht]
\begin{center}
\includegraphics[width=120mm]{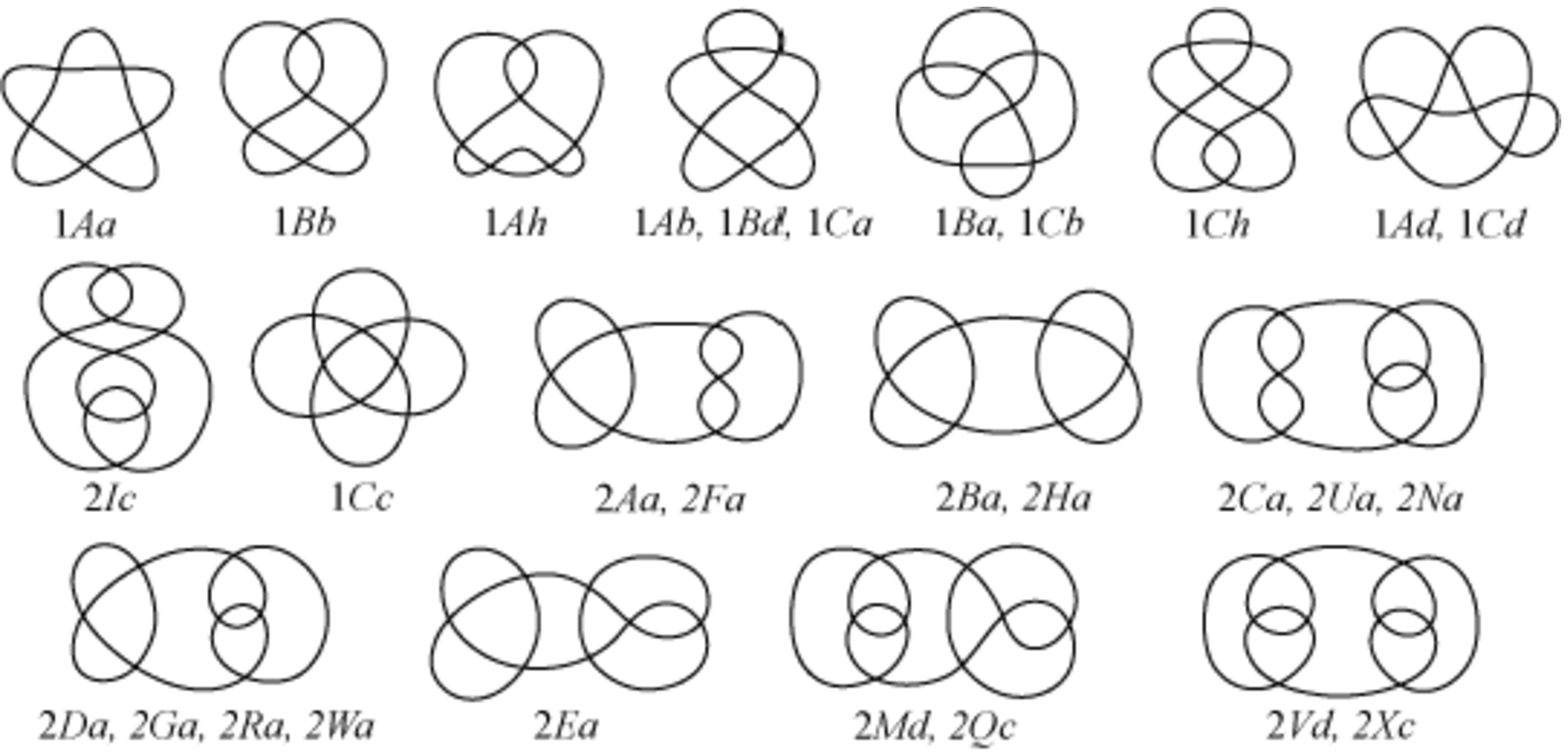}
\caption{All the reduced knot projections with warping degree two. }
\label{mwd2kp}
\end{center}
\end{figure}

\noindent Now, we have obtained all the prime alternating knots of minimal warping degree one and two, and we have the following. 

\phantom{x}
\begin{corollary}
We have the following: \\
\noindent (1) $md(3_1)=md(4_1)=1$ \\
\noindent (2) $md(5_1)=md(5_2)=2$ \\
\noindent (3) $md(6_1)=md(6_2)=md(6_3)=2$\\
\noindent (4) $md(7_i)=3$ for $1 \leq i \leq 5$ \\
\noindent (5) $md(7_6)=md(7_7)=2$ \\
\noindent (6) $md(8_j)=3$ for $1 \leq j \leq 11$ or $13 \leq j \leq 17$\\
\noindent (7) $md(8_{12})=md(8_{18})=2$ \\
\noindent (8) $md(K) \geq 3$ for any prime alternating knot $K$ with $c(K) \geq 9$
\end{corollary}
\phantom{x}

\begin{proof}
(1) is shown in \cite{JS}. 
(2), (3), (5), (7) and (8) follow from Theorem \ref{prime-alt-two}. 
(4) and (6): Let $D$ be an oriented minimal diagram of a prime alternating knot $K$. 
Since $D$ is alternating, $d(D)+d(-D)=c(D)-1$ holds (Theorem 1.3 in \cite{shimizu-wd}), i.e., $d(D)+d(-D)=c(K)-1$ holds. 
If $c(K)=7$ and $md(K)>2$, we have $d(D)=d(-D)=3$ from $d(D)+d(-D)=6$ and $d(D), d(-D) \geq 3$. 
If $c(K)=8$ and $md(K)>2$, we have $d(D)=3$ and $d(-D)=4$ or $d(D)=4$ and $d(-D)=3$ from $d(D)+d(-D)=7$ and $d(D), d(-D) \geq 3$. 
\end{proof}

\section{Estimation of the warping degree of a knot projection by r-factors}

In this section, we define and discuss the ``r-factor'' to improve Corollary \ref{k1k2k3p} and give an estimation for the warping degree of a knot projection. 
Let $P$ be a knot projection, and let $c$ be a crossing of $P$. 
We denote the pair of $P$ and $c$ by $P^c$. 
A set $S$ of regions of $P$ is an {\it r-factor of $P^c$} when $S$ satisfies the following three conditions: 
(i) each crossing of $P$ except $c$ is on the boundary of one of the regions in $S$, 
(ii) the crossing $c$ is not on the boundary of any region in $S$, 
(iii) all the regions in $S$ are disjoint. 
Some examples are shown in Figure \ref{52}. 
\begin{figure}[ht]
\begin{center}
\includegraphics[width=80mm]{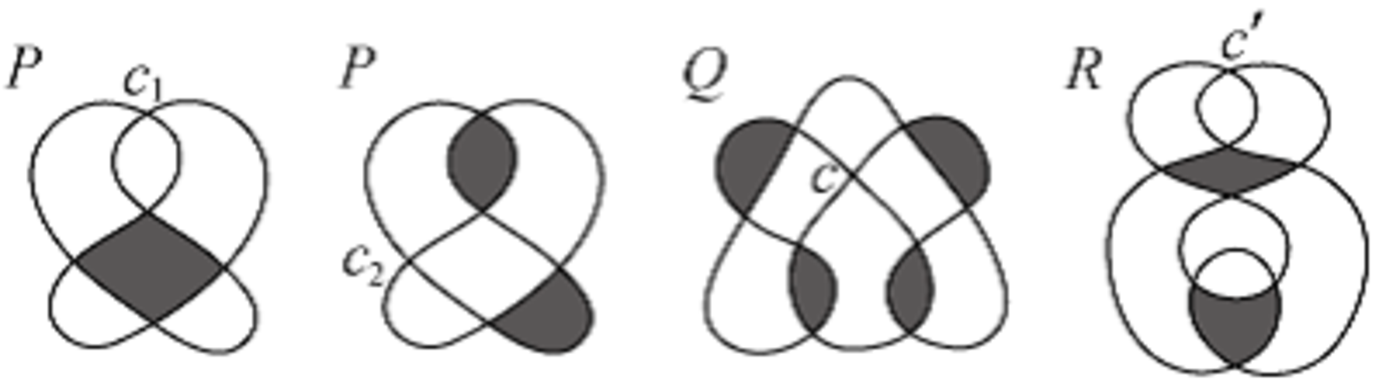}
\caption{The r-factors of $P^{c_1}, P^{c_2}, Q^c$ and $R^{c'}$.  }
\label{52}
\end{center}
\end{figure}
We note that some knot projections have no r-factor with any crossing. 
For example, see the knot projection $P$ in Figure \ref{818}. 
Take any crossing $c$ without loss of generality by the symmetry of the knot projection $P$. 
We need the region $R_1$ or $R_2$ to cover the crossing $c_1$. 
If we choose $R_1$, we can not choose $R_2, R_4$ and $R_5$. 
If we choose $R_2$, we can not choose $R_1, R_2, R_4, R_5$ and $R_6$. 
Thus we can not cover all crossings. 
Hence $P$ has no r-factors. 
In this section, we show the following. 
\begin{figure}[ht]
\begin{center}
\includegraphics[width=25mm]{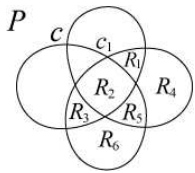}
\caption{The knot projection $P$ has no r-factors with any crossing. }
\label{818}
\end{center}
\end{figure}

\phantom{x}
\begin{theorem}
Let $P$ be a knot projection without 1-gons, and let $c$ be a crossing of $P$. 
If a set of $n$ regions of $P$ consisting of $k_1$-gon, $k_2$-gon, $\dots$ ,$k_n$-gon is an r-factor of $P^c$, the following inequality holds: 
$$n \leq d(P) \leq k_1 + k_2 + \dots + k_n -n.$$
\label{ndPk}
\end{theorem}

\noindent For example, we have $1 \leq d(P) \leq 3$ from the r-factor of $P^{c_1}$ in Figure \ref{52}. 
We have $2 \leq d(R) \leq 5$ from the r-factor of $R^{c'}$ in Figure \ref{52}.\footnote{An alternating diagram obtained from $R$ has warping degree 2 and 5 for each orientation. } 
In particular, for the case that $k_1 = k_2 = \dots =k_n =2$, we have the following. 

\phantom{x}
\begin{corollary}
Let $P$ be a knot projection without 1-gons. 
If $P$ has an r-factor of $P^c$ which consists of $n$ 2-gons for a crossing $c$, then $d(P)=n$. 
\label{ndPk-c}
\end{corollary}
\phantom{x}

\noindent For example, $P^{c_2}$ and $Q^c$ in Figure \ref{52} has an r-factor consisting of 2-gons, and we obtain $d(P)=2$ and $d(Q)=4$. 
To prove Theorem \ref{ndPk}, we show the following lemma. 

\phantom{x}
\begin{lemma}
Let $D$ be an oriented alternating knot diagram which has a $k$-gon and a crossing $c$ which is not on the $k$-gon. 
If we take a base point $b$ just before an over-crossing of $c$, then the $k$-gon has at most $k-1$ warping crossing points of $D_b$ on the boundary. 
\label{upper-wd}
\end{lemma}
\phantom{x}

\begin{proof}
Let $T$ be a $k$-tangle which bounds a regular neighborhood of the $k$-gon. 
Since $T$ has $k$ crossings, it is obvious that $T$ has at most $k$ warping crossing points of $D_b$. 
Let $S$ be the segment of $T$ such that we encounter $S$ first from $b$. 
Two crossings belong to $S$, one over-crossing and one under-crossing. 
The over-crossing is not a warping crossing point of $D_b$. 
Hence $T$ has at most $k-1$ warping crossing points of $D_b$. 
\end{proof}
\phantom{x}

\noindent We show Theorem \ref{ndPk}. 

\phantom{x}
\noindent {\it Proof of Theorem \ref{ndPk}.} \ 
Let $D$ be an oriented alternating knot diagram with $|D|=P$ and $d(D)=d(P)$. 
From Corollary \ref{k1k2k3p} and Lemma \ref{upper-wd}, we have $n \leq d(D) \leq (k_1-1)+(k_2-2)+\dots +(k_n-1)$. 
\hfill$\Box$

\phantom{x} 
\noindent From Corollary \ref{ndPk-c}, we have the following corollary. 

\phantom{x}
\begin{corollary}
(1) Let $P$ be a knot projection represented by the Conway notation with a positive odd integer $k$. 
Then $d(P)= (k-1)/2$. \\
\noindent (2) Let $P$ be a knot projection represented by the Conway notation with two integers $lm$, where $l$ is a positive even integer and $m$ is a positive odd integer. 
Then $d(P)= l/2 +(m-1)/2$. 
\label{md-conway}
\end{corollary}
\phantom{x}

\begin{proof}
For both cases (1) and (2), $P$ has an r-factor consisting of 2-gons with a crossing. 
\end{proof}
\phantom{x}

\noindent From Corollary \ref{md-conway}, we obtain $md(K)=(k-1)/2$ for a $(2,k)$-torus knot $K$ and $md(L)=(m+1)/2$ for a twist knot $L$ with $c(L)=m+2$. 
Remark that each $(2,k)$-torus knot and twist knot has a unique minimal diagram (see, for example, \cite{shimizu-rcc}). 

\section*{Acknowledgment}
The author is deeply grateful to Slavik Jablan for valuable discussions on the minimal warping degree.


\begin{thebibliography}{99}
\bibitem{JS}S.~Jablan and A.~Shimizu, On the warping sum of knots, J. Knot Theory Ramifications {\bf 27} (2018), 1843002 [8 pages]. 
\bibitem{kauffman} L.~Kauffman, State models and the Jones polynomial, Topology {\bf 26} (1987), 395--407. 
\bibitem{kawauchi-lecture} A.~Kawauchi, Lectures on knot theory (in Japanese), Kyoritsu shuppan Co. Ltd, 2007. 
\bibitem{LM} W.~B.~R.~Lickorish and K.~C.~Millett, A polynomial invariant of oriented links, Topology {\bf 26} (1987), 107--141. 
\bibitem{murasugi} K.~Murasugi, Jones polynomials and classical conjectures in knot theory, Topology {\bf 26} (1987), 187--194. 
\bibitem{murasugi-book} K.~Murasugi, Knot theory and its applications (Birkh\"{a}user, 1996). 
\bibitem{MO} M.~Ozawa, Ascending number of knots and links, J. Knot Theory Ramifications {\bf 9} (2010), 15--25. 
\bibitem{shimizu-rcc}A.~Shimizu, Region crossing change is an unknotting operation, J. Math. Soc. Japan {\bf 66} (2014), 693--708. 
\bibitem{shimizu-wd}A.~Shimizu, The warping degree of a knot diagram, J. Knot Theory Ramifications {\bf 19} (2010), 849--857. 
\bibitem{shimizu-w-poly}A.~Shimizu, The warping polynomial of a knot diagram, J. Knot Theory Ramifications {\bf 21} (2012), 1250124. 
\bibitem{thistlethwaite} M.~Thistlethwaite, A spanning tree expansion of the Jones polynomial, Topology {\bf 26} (1987), 297--309. 
\end{thebibliography}
\end{document}